\documentclass[12pt]{article}
\usepackage{amsmath, amssymb, graphics,graphicx,geometry}
\newtheorem{theorem}{Theorem}

\newtheorem{definition}{Definition}

\newtheorem{lemma}{Lemma}

\title{Proving the existence of the $n$th root by induction.}

\author{Alvaro Salas\thanks{Department of Mathematics, Universidad de
Caldas, Department of Mathematics, Universidad Nacional de Colombia,
Manizales, Caldas. \emph{email} : asalash2002@yahoo.com}}
\date{}
\begin {document}
\maketitle
\begin {abstract}
In this paper we prove by induction on $n$ that any positive real
number has $n$th root.
\end{abstract}
\section{Auxiliary facts}
We define the set $\mathbb{R}$ of real numbers as a numeric ordered
field in which the following axiom holds:\\
\textbf{Least Upper Bound Axiom }: If $\emptyset\neq E\subseteq
\mathbb{R}$ is a non empty set bounded from above, then there exists
a least upper bound for $E$. \\The least upper bound of $E$ is
unique and it is denoted by $\sup E$. If $x=\sup E$, then :\\
\textbf{A}. $t\leq x$ for any $t\in E$, that is, $x$ is an upper bound for $E$.\\
\textbf{B}. For any $\delta>0$ we may find $t\in E$ such that
$x-\delta<t\leq x$.
\begin{definition}
Let $n$ be a positive integer and $a>0$. We say that $x$ is a $n$th
root of $a$ if $x^n=a$.
\end{definition}
It is easy to show that if $a$ has an $n$th root, then this root is
unique. This follows from the fact that if $x$ and $y$ are positive
numbers for which $x^n=y^n$, then $x=y$. The $n$th root of $a$ is
denoted by $\sqrt[n]{a}$.
\begin{lemma}\label{lema0} Suppose that  $a>0$ and $n$ is a positive integer. The set
$$E=\{t\geq 0\,\,|\,t^n<a\}.$$ is not empty and bounded from above.
\end{lemma}
\textbf{Proof}. $E$ is not empty, because $0\in E$. On the other
hand, it is easy to show by induction on $n$ that $(a+1)^n>a$. Now,
if $t>a+1$ then $t^n>(a+1)^n>a$ and then $t\notin E.$ This means
that $a+1$ is an upper bound for $E$.
\section{Existence of $n$th root of a positive number}
\begin{theorem}\label{teo1}
For any positive integer $n\geq 1$ and for any $a>0$ there exists
$\sqrt[n]{a}>0.$
\end{theorem}
 We first prove by induction  that for any
positive integer $n$ and for any $b>0$ there exists $\sqrt[2^n]{b}.$
This is true for $n=1$. Indeed, according to Lemma \ref{lema0}, the
set $E=\{t>0\,\,|\,t^2<a\}$ has a least upper bound $x=\sup E\geq
0$. We claim that $x^2=a$. Indeed, let $0<\varepsilon<1$ and define
$\delta=\varepsilon/(2x+1).$ Clearly, $0<\delta\leq \varepsilon<1$.
There exists an element $t\in E$ such that $x-\delta<t\leq x$. We
have
$$x^2<(t+\delta)^2=t^2+2t\delta+\delta^2<a+(2t+1)\delta\leq a+\varepsilon.$$
Since the inequality $x^2<a+\varepsilon$ holds for any
$\varepsilon\in(0,1)$, we conclude that $x^2\leq a$. \\
On the other hand, $x+\delta\notin E$, so
$$a\leq (x+\delta)^2<x^2+(2x+1)\delta=x^2+\varepsilon.$$
Since the inequality $x^2+\varepsilon>a$ holds for any
$\varepsilon\in(0,1)$, we conclude that $x^2\geq a$. We have proved
the equality $x^2=a$. In particular, $x=\sqrt{a}>0$.\\
Suppose that for some positive integer $k$ the number
$\sqrt[2^k]{c}$ is defined for every $c>0$. Let $b$ be a positive
number. By inductive hypothesis, for $c=\sqrt{b}$ there exists $y>0$
such that $y^{2^{k}}=\sqrt{b}$. From this it follows that
$$b=\sqrt{b}^2=y^{2\cdot2^{k}}=y^{2^{k+1}}$$
and then $y=\sqrt[2^{k+1}]{b}$, that is, our Lemma is also true for
$n=k+1$.\\
Finally, we shall prove that for any positive integer $n\geq 2$ and
for any $a>0$ there exists $\sqrt[n]{a}>0.$ We proceed by induction
on $n$. Our theorem is true for $n=2$ since any positive number a
has a square root. Suppose that $\sqrt[n-1]{c}$ exists for any $c>0$
and for some $n\geq 3$. Let $a>0$ and define the set
$$E=\{t\geq 0\,\,|\,t^n<a\}.$$
By Lemma \ref{lema0} this set has a supremum, say $x=\sup E\geq 0$.
 Define $m=2^n-n$. Then $m+n=2^n$ and $m$ is a
positive integer, since $2^n>n$. There exists a positive $y$ for
which $y^{2^n}=y^{m+n}=ax^m$. Observe that
\begin{equation}\label{ueq}
y^n=y^{-m}y^{m+n}=y^{-m}ax^m=\left(\dfrac{x}{y}\right)^ma
\end{equation}
We claim that $x^n=a$.\\
Indeed, suppose that $x^n<a$. Then $x^{m+n}<ax^m=y^{m+n}$. This
implies that $x<y$. On the other hand, from (\ref{ueq}), $y^n<a$ and
then $y\in E$, so $y\leq x$. But $y>x$ and we get a
contradiction.\\
Now, suppose that $x^n>a$. Then $x^{m+n}>ax^m=y^{m+n}$ and this
implies that $x>y$. On the other hand, (\ref{ueq}) implies that
$y^n>a$.  Now, if $t\in E$, then $t^n<a<y^n$. This inequality
implies that $t<y$ for any $t\in E$. This says us that $y$ is an
upper bound of $E$ that is \emph{less} than $x=\sup E$ and we again
get a
contradiction.\\
We have proved the equality $x^n=a$. In  particular
$x=\sqrt[n]{a}>0$.
\end{document}